\documentclass[11pt]{article}

\usepackage[utf8]{inputenc}
\usepackage{amsmath}
\usepackage{amsfonts,amssymb,latexsym,multicol,color}
\usepackage[francais]{babel}
\usepackage[T1]{fontenc}

\newcounter{fig}
\newtheorem{theo}{Th\'eor\`eme}

\newtheorem{prop}{Proposition}

\makeatletter
\ifcase\@ptsize

\or

\or

\fi
\makeatother

\newcommand{\fhi}{\varphi}
\newcommand{\ioe}{\leqslant}
\newcommand{\soe}{\geqslant}

\newcommand{\vers}{\rightarrow}

\newcommand{\demi}{{\frac{1}{2}}}

\newcommand{\Bcal}{{\mathcal B}}

\newcommand{\Nat}{{\mathbb N}}

\newcommand{\Rat}{{\mathbb Q}}

\newcommand{\fin}{\hfill$\Box$}
\newcommand{\dem}{\noindent {\bf D\'emonstration\ }}

\providecommand{\bysame}{\leavevmode ---\ }
\providecommand{\og}{``}
\providecommand{\fg}{''}
\providecommand{\smfandname}{et}

\newcommand{\virg}{\raisebox{.7mm}{,}}

\title{Sur le minimum de la fonction de Brjuno}

\author{Michel Balazard et Bruno Martin}

\begin{document}
\maketitle

\begin{center}
  {\sc Abstract}
\end{center}
\begin{quote}
{\footnotesize The Brjuno function attains a strict global minimum at the golden section.}
\end{quote}

\begin{center}
  {\sc Keywords}
\end{center}
\begin{quote}
{\footnotesize Gauss transformation, Brjuno function, Golden section }
\end{quote}



La transformation de Gauss,
\[
\alpha(x)=\{1/x\},
\]
où $\{t\}$ désigne la partie fractionnaire du nombre réel $t$, est bien définie sur l'ensemble
\[
X=\, ]0,1[\,\setminus\, \Rat,
\]
et à valeurs dans $X$. On peut donc considérer ses itérées successives, définies par les relations~$\alpha_0(x)=x$ et $\alpha_{k+1}(x)=\alpha_k\big(\alpha(x)\big)$ pour $k \soe 0$, et les produits
\[
\beta_k(x)=\prod_{j=0}^k \alpha_j(x) \quad (k \in \Nat, \, x \in X),
\]
avec la convention supplémentaire $\beta_{-1}(x)=1$. Comme $\alpha$ est continue sur $X$, les $\alpha_k$ et les $\beta_k$ le sont également. 

\smallskip

La fonction de Brjuno est alors définie, pour tout $x \in X$, comme la somme, éventuellement égale à $+\infty$, de la série à termes positifs
\[
\Phi(x)=\sum_{k \soe 0}\beta_{k-1}(x)\ln\big(1/\alpha_k(x)\big).
\]
Les points de convergence sont appelés \emph{nombres de Brjuno} ; nous noterons $\Bcal$ leur ensemble. La fonction et les nombres de Brjuno interviennent dans la théorie des systèmes dynamiques (cf. par exemple \cite{buff-cheritat-2006}, \cite{marmi-moussa-yoccoz-1997}, \cite{zbMATH05183581}, \cite{zbMATH00827357}). L'ensemble $\Bcal$ est de mesure $1$, et donc dense dans $[0,1]$.

L'ensemble $\Bcal$ des nombres de Brjuno est stable par $\alpha$. La fonction $\Phi$ vérifie l'équation fonctionnelle
\begin{equation}
  \label{181029a}
\Phi(x)=\ln(1/x) +x \Phi\big(\alpha(x)\big) \quad (x \in \Bcal),  
\end{equation}
et, plus généralement,
\begin{equation}\label{181027a}
\Phi(x)=\Phi_K(x)+\beta_K(x)\Phi\big(\alpha_{K+1}(x)\big) \quad (K \in \Nat, \,  x \in \Bcal),
\end{equation}
où $\Phi_K$ désigne la somme partielle
\[
\Phi_K(x)=\sum_{k=0}^K\beta_{k-1}(x)\ln\big(1/\alpha_k(x)\big).
\]
Les fonctions $\Phi_K$ sont définies et continues sur $X$.

\smallskip

Dans l'article \cite{rivoal-2010}, Rivoal émet plusieurs conjectures sur les valeurs extrémales de séries diophantiennes, dont certaines sont proches de la fonction de Brjuno. Le théorème suivant fournit la réponse à une question posée aux auteurs par Rivoal, concernant la fonction~$\Phi$ elle-même.

\begin{theo}
Soit $\theta=(\sqrt{5}-1)/2=0,618\dots$ le nombre d'or. Pour tout nombre de Brjuno $x\neq \theta$, on a $\Phi(x) > \Phi(\theta)$.
\end{theo}



\medskip 

La démonstration de ce théorème s'appuie sur cinq propositions auxiliaires.

\begin{prop}\label{181027b}
Soit $r$ un nombre rationnel, élément du segment $[0,1]$. On a alors
\[
\Phi(x) \vers \infty \quad (x\vers r, \, x \in \Bcal).
\]
\end{prop}
\dem

On a $\Phi(x) \soe \Phi_0(x)=\ln 1/x$, donc le résultat est vrai si $r=0$. On a
\[
\Phi(x)=\ln 1/x + x \Phi\big(\alpha(x)\big) \soe \demi\Phi(1/x-1) \quad ( 1/2 < x <1),
\]
donc le résultat est aussi vrai si $r=1$.

Si $r$ est un nombre rationnel de $]0,1[$, écrivons $r$ sous forme d'une fraction continue finie,
\[
r=[0;a_1,\dots,a_k] \quad (k\soe 1, \, a_1, \dots, a_k \in \Nat^*, \, a_k \soe 2).
\]
L'application $\fhi : \,t \mapsto [0;a_1,\dots,a_{k-1},a_k+t]$ est un homéomorphisme de $]-1,1[$ sur un certain voisinage de $r$ dans $]0,1[$. 

Pour $t\in X$, ce qui entraîne $\fhi(\pm t) \in X$, on a
\[
\alpha_k\big(\fhi(t)\big)=t \quad ; \quad \alpha_k\big(\fhi(-t)\big)=1-t,
\]
et $\beta_{k-1}\big(\fhi(\pm t)\big)$ est minorée par une constante positive (dépendant de $r$, cf. \cite{Brjuno}, (18)-(19), p. 199) pour $t \in X$. Comme
\[
\Phi\big(\fhi(t)\big)\soe \beta_{k-1}\big(\fhi( t)\big)\Phi\Big(\alpha_{k}\big(\fhi(t)\big)\Big) \quad ( t  \in\, ]-1,1[\,\setminus\Rat),
\]
le résultat général découle des cas particuliers $r=0,1$.\fin

\smallskip

Posons maintenant
\[
C=\inf_{x \in \Bcal} \Phi(x).
\]
 
On a $0 \ioe C < \infty$, car $\Phi$ est à valeurs $\soe 0$ et $\Bcal$ n'est pas vide.

\begin{prop}\label{181027f}
La borne inférieure $C$ est le minimum de la fonction de Brjuno sur $\Bcal$ : il existe~$r \in \Bcal$ tel que $C=\Phi(r)$. \end{prop}
\dem

Soit $(x_n)_{n\soe 1}$ une suite d'éléments de $\Bcal$ telle que
\[
\Phi(x_n) \vers C \quad ( n \vers \infty).
\]
En remplaçant éventuellement cette suite par une de ses sous-suites, nous pouvons supposer que la suite $(x_n)$ est elle-même convergente, vers une limite $r \in [0,1]$. La proposition \ref{181027b} et l'hypothèse de convergence vers $C$ de la suite $\Phi(x_n)$ entraînent alors que $r$ est irrationnel. Nous allons montrer que $r \in \Bcal$ et que $\Phi(r)=C$.

Soit $K \in \Nat$. Par définition de $\Phi$ on a $\Phi(x_n) \soe \Phi_K(x_n)$ pour tout $n$. La continuité de $\Phi_K$ en tout point irrationnel entraîne donc, par passage à la limite, l'inégalité $C \soe \Phi_K(r)$. Comme $K$ est arbitraire, cela prouve que $r$ est un nombre de Brjuno, et que $C \soe \Phi(r)$. Comme on a aussi, par définition de $C$, l'inégalité inverse $\Phi(r) \soe C$, on en déduit l'égalité $C=\Phi(r)$. Cette borne inférieure est donc bien un minimum.\fin

\begin{prop}\label{181027c}
Soit $r \in \Bcal$ tel que $C=\Phi(r)$. Pour tout $K \in \Nat$, on a
\[
C=\Phi(r) \soe \frac{\Phi_K(r)}{1-\beta_K(r)}\cdotp
\]
\end{prop}
\dem

En appliquant la relation \eqref{181027a} à $x=r$, nous obtenons
\[
C=\Phi(r)=\Phi_K(r)+\beta_K(r)\Phi\big(\alpha_{K+1}(r)\big)\soe \Phi_K(r)+C\beta_K(r),
\]
par définition de $C$. L'assertion en résulte.\fin

\begin{prop}\label{181024b}
Pour tout $K \in \Nat$, on a
\begin{equation*}
\Phi(\theta)=\frac{\Phi_K(\theta)}{1-\beta_K(\theta)}\cdotp
\end{equation*}
\end{prop}
\dem

En effet, en appliquant la relation \eqref{181027a} à $x=\theta$, nous obtenons
\begin{equation*}
\Phi(\theta)=\Phi_K(\theta) +\beta_K(\theta)\Phi\big(\alpha_{K+1}(\theta)\big)=\Phi_K(\theta) +\beta_K(\theta)\Phi(\theta),
\end{equation*}
puisque $\theta$ est point fixe de $\alpha$.\fin

\begin{prop}\label{181027e}
Soit $r \in \Bcal$ tel que $C=\Phi(r)$. On a $r \soe \theta$.
\end{prop}
\dem

D'après les propositions \ref{181027c} et \ref{181024b} avec $K=0$, et la définition de $C$, on a
\[
\frac{\Phi_0(\theta)}{1-\beta_0(\theta)}=\Phi(\theta) \soe C=\Phi(r) \soe \frac{\Phi_0(r)}{1-\beta_0(r)}\cdotp
\]

Or, pour $0<x<1$, 
\[
\frac{\Phi_0(x)}{1-\beta_0(x)}=\frac{\ln 1/x}{1-x}=\int_0^1\frac{dt}{(1-t)x +t}\virg
\]
est une fonction strictement décroissante de $x$, donc $\theta \ioe r$.\fin

\medskip

Au vu de la proposition \ref{181027f}, l'énoncé suivant est équivalent à celui du théorème.

\begin{prop}\label{191105e}
Soit $r \in \Bcal$ tel que $C=\Phi(r)$. On a $r = \theta$.
\end{prop}
\dem

D'après les propositions \ref{181027c} et \ref{181024b} avec $K=1$, et la définition de $C$, on a
\begin{equation}\label{181027d}
\frac{\Phi_1(\theta)}{1-\beta_1(\theta)}=\Phi(\theta) \soe C=\Phi(r) \soe \frac{\Phi_1(r)}{1-\beta_1(r)}\cdotp
\end{equation}

Pour $1/2<x<1$, posons
\[
f(x)=\frac{\Phi_1(x)}{1-\beta_1(x)}=\frac{\ln 1/x +x \ln1/\alpha(x)}{1-(1-x)}=\frac{\ln 1/x}x+\ln \frac{x}{1-x}\virg
\]
puisque $\alpha(x)=(1-x)/x$ pour $1/2 < x <1$.

On a
\[
f'(x)=\frac{\ln x}{x^2}+ \frac{2x-1}{x^2(1-x)} \virg
\]
fonction qui a le signe de
\[
g(x)=2x-1+(1-x)\ln x.
\]

La fonction $g$ est strictement croissante sur $]0,1]$ et $g(0,61)=0,027\!\ldots>0$. Par conséquent, la fonction $f$ est strictement croissante sur $[\theta, 1[$. La relation \eqref{181027d} et la proposition \ref{181027e} entraînent donc l'égalité $r=\theta$.\fin

\medskip

Pour mettre en évidence la substance de la démonstration que nous venons d'exposer, remplaçons dans la définition de la fonction de Brjuno, la fonction logarithme par une fonction arbitraire $u :\, [1,\infty[ \vers [0,\infty[$, et posons
\[
\Psi_u(x)=\sum_{k \soe 0}\beta_{k-1}(x)\,u\big(1/\alpha_k(x)\big)
\]
(cf. \cite{zbMATH05183581}, \S3.2, p. 609). Comme tout irrationnel quadratique appartient à l'ensemble $\Bcal_u$ des points de convergence de cette série\footnote{D'une part, on a l'inégalité $\beta_{k-1} \ioe 1/F_{k+1}$, où $F_n$ désigne le $n$\up{e} nombre de Fibonacci ; d'autre part, lorsque~$x$ est un irrationnel quadratique, l'ensemble des valeurs $\alpha_k(x)$ ($k \in \Nat$), est fini.}, cet ensemble est toujours dense. L'équation fonctionnelle \eqref{181029a} devient
\begin{equation}
  \label{191106a}
\Psi_u(x)=u(1/x) +x \Psi_u\big(\alpha(x)\big) \quad (x \in \Bcal_u),  
\end{equation}

En toute généralité, il n'est pas vrai que la fonction~$\Psi_u$ atteigne un minimum global au nombre d'or. En effet, si $u(t)=t^a$, avec $a <1$, et $\theta'=1-\theta=1/(2+\theta)$, de sorte que~$\alpha(\theta')=\theta$, on a
\[
\Psi_u(\theta)=\theta^{-a}+\theta\Psi_u(\theta)=\theta^{-a}+(1-\theta')\Psi_u(\theta), 
\]
donc $\theta'\Psi_u(\theta)=\theta^{-a}$, et
\[
 \Psi_u(\theta')=\theta'^{-a}+\theta'\Psi_u(\theta)=\theta'^{-a}+\theta^{-a}<\theta^{-a}/\theta'=\Psi_u(\theta),
 \]
car la dernière inégalité équivaut à $\theta'^{1-a} < \theta^{1-a}$, c'est-à-dire à $\theta'<\theta$, inégalité vraie.

\medskip

En remplaçant, dans les énoncés, $\Phi$ par $\Psi_u$, $\Bcal$ par $\Bcal_u$, et $\Phi_K$ par la somme partielle de $\Psi_u$ correspondante, examinons les conditions de validité des six propositions précédentes.

\smallskip

$\bullet$ La proposition \ref{181027b} reste valable si, et seulement si 
\begin{equation}\label{191105a}
u(x) \vers \infty \quad (x \vers \infty).
\end{equation}

$\bullet$ Sous l'hypothèse \eqref{191105a}, la proposition \ref{181027f} reste valable si
\begin{equation}\label{191105b}
\text{la fonction $u$ est continue en chaque point irrationnel.} 
\end{equation}

$\bullet$ Sous les hypothèses \eqref{191105a} et \eqref{191105b}, les propositions \ref{181027c} et \ref{181024b} restent valables.

$\bullet$ Sous les hypothèses \eqref{191105a} et \eqref{191105b}, la proposition \ref{181027e} reste valable si
\begin{equation}\label{191105c}
\text{la fonction } x \mapsto \frac{u(1/x)}{1-x} \text{ est strictement décroissante sur } ]0,1[.
\end{equation} 

$\bullet$ Sous les hypothèses \eqref{191105a}, \eqref{191105b} et \eqref{191105c}, la proposition \ref{191105e} reste valable si
\begin{equation}\label{191105d}
\text{la fonction } x \mapsto \frac{u(1/x)}{x}+u\Big(\frac{x}{1-x}\Big) \text{ est strictement croissante sur } [\theta,1[.
\end{equation} 

\smallskip

Notons que l'hypothèse \eqref{191105c} entraîne que $u(x)$ tend vers $0$ quand $x$ tend vers $1$. La proposition suivante décrit une classe de fonctions (comprenant la fonction logarithme) vérifiant les hypothèses \eqref{191105a}, \eqref{191105b}, \eqref{191105c} et \eqref{191105d} ; les fonctions~$\Psi_u$ correspondantes ont donc toutes un minimum  global strict au nombre d'or.

\begin{prop}\label{191106b}
Soit $u :[1,\infty[\, \vers [0,\infty[$ une fonction telle que 

$\bullet$ $u(1)=0$ ;

$\bullet$ $u(x) \vers \infty \quad (x \vers \infty)$ ;

$\bullet$ la fonction~$v$ définie par $v(t)= tu(t)$ est strictement convexe sur $[1,\infty[$.\\
Alors les hypothèses \eqref{191105a}, \eqref{191105b} et \eqref{191105c} sont vérifiées. Si, de plus, 

$\bullet$ la fonction $x \mapsto u\big(x/(1-x)\big)$ est convexe sur $[1/2,1[$ ;

$\bullet$ $u'(1/\theta) \soe u(1/\theta)$\\
(où $u'$ désigne la dérivée à droite), alors l'hypothèse \eqref{191105d} est vérifiée.
\end{prop}
\dem

La convexité de $v$ entraîne sa continuité, donc aussi celle de $u$, sur $]1,\infty[$. Les deux conditions~$u\soe 0$ et $u(1)=0$  entraînent aussi la continuité à droite au point $1$.

La stricte convexité de $v$ et le fait que $v(1)=0$ entraînent que la pente $w(t)=v(t)/(t-1)$ est une fonction strictement croissante de $t>1$. Par conséquent la fonction $v$ est elle-même strictement croissante, ainsi que la fonction $u(t)=(1-1/t)w(t)$. De plus, la fonction composée
\[
x \mapsto w(1/x)=\frac{u(1/x)}{1-x}
\]
est strictement décroissante sur $]0,1[$.  

La stricte croissance et la convexité de $v$ sur $[1,\infty[$, et la stricte convexité de la fonction définie par $x \mapsto 1/x$ sur $]0,1[$, entraînent la stricte convexité de la fonction $x\mapsto v(1/x)=u(1/x)/x$ sur~$]0,1[$.

Si la fonction composée $x \mapsto u\big(x/(1-x)\big)$ est convexe sur son domaine de définition $[1/2,1[$, alors la fonction intervenant dans \eqref{191105d} est somme d'une fonction strictement convexe et d'une fonction convexe sur $[1/2,1[$ ; elle est donc strictement convexe sur cet intervalle. Pour vérifier~\eqref{191105d}, il suffit de vérifier que sa dérivée à droite au nombre d'or est positive ou nulle. Cela s'écrit
\[
-\frac{1}{\theta^2}u(1/\theta) -\frac{1}{\theta^3}u'(1/\theta)+\frac{1}{(1-\theta)^2}u'\Big(\frac{\theta}{1-\theta}\Big) \soe 0
\]

En tenant compte de l'équation $\theta/(1-\theta)=1/\theta$, on voit que cette inégalité équivaut à celle de l'énoncé.\fin

\smallskip

Si $u : [1,\infty[ \, \vers [0,\infty[$ est strictement convexe, vérifie $u(1)=0$ et $u'(1/\theta) \soe u(1/\theta)$, alors elle vérifie les hypothèses de la proposition \ref{191106b}. En particulier, c'est le cas des fonctions $t \mapsto (t-1)^a$, avec $a > 1$ (si $a=1$, la proposition \ref{191106b} s'applique directement). Par ailleurs, les fonctions $t \mapsto \ln^at$, avec $a\soe 1$, sont des exemples de fonctions non convexes vérifiant également cette proposition.


\medskip
\begin{center}
{\sc Remerciements}
\end{center}

{\footnotesize Nous remercions Tanguy Rivoal d'avoir suscité le présent travail, et l'arbitre anonyme d'avoir suggéré de généraliser notre étude. Cette recherche a été rendue possible par le programme \emph{Research in pairs} du Mathematisches Forschungsinstitut Oberwolfach. Nous remercions cette institution pour les conditions de travail idéales dont nous avons bénéficié.}


\medskip

\footnotesize

\noindent BALAZARD, Michel\\
Aix Marseille Univ, CNRS, Centrale Marseille, I2M, Marseille, France\\
Adresse \'electronique : \texttt{balazard@math.cnrs.fr}

\medskip

\noindent MARTIN, Bruno\\
ULCO,  LMPA, Calais, France\\
Adresse \'electronique : \texttt{Bruno.Martin@univ-littoral.fr}

\end{document}